\documentclass[12pt,a4paper,twoside,reqno]{amsart}

\usepackage[english]{babel}

\usepackage{amsmath}             \usepackage{amssymb}
\usepackage{amsthm}              \usepackage{latexsym}
\usepackage{textcomp}            \usepackage{graphicx}

\newtheorem{thm}{Theorem}        \newtheorem{rem}{Remark}
\newtheorem{lem}{Lemma}          \newtheorem{prop}{Proposition}
\theoremstyle{definition}        \newtheorem{Def}{Definition}

\newcommand{\Ro}{{\mathbb{R}}}   \newcommand{\Co}{{\mathbb{C}}}
\newcommand{\eps}{{\varepsilon}}

\title[Numb. of Lim. Cycl. of Even Deg. Li\'enard Eqns for Foci]
{One Upper Estimate\\
on the Number of Limit Cycles\\
of Even Degree Li\'enard Equations\\
in the Focus Case}

\author[Grisha Kolutsky]{Grisha Kolutsky}

\date{18 November 2009}

\thanks{The work exposed here was partially supported by
grants {7-01-00017-a} and {08-01-00342-a} of the Russian
Foundation for Basic Research, by the grant\linebreak
{No.~NSh-3038.2008.1} of the President of Russia for support of
leading scientific schools and by the Russian Universities grant
{No.~RNP.2.1.1.5055.}}

\address{\rm Grisha Kolutsky\newline \it
Department of the Theory of Dynamical Systems\newline
Faculty of Mechanics and Mathematics\newline
Lomonosov Moscow State University\newline
MSU, GSP, Glavnoe Zdanie, Leninskie Gory\newline
119899 Moscow, Russia\newline
{\rm e-mail: kolutsky AT mccme DOT ru}}

\begin{document}

\begin{abstract}
We give an explicit upper bound for a number of limit cycles of
the Li\'enard equation $\dot{x}=y-F(x)$, $\dot{y}=-x$ of even
degree in the case its unique singular point $(0,0)$ is a focus.

M.~Caubergh and F.~Dumortier get an explicit linear upper estimate
for the number of large amplitude limit cycles of such equations
\cite{CD}. We estimate the number of mid amplitude limit cycles of
Li\'enard equations using the Growth-and-Zeros theorem proved by
Ilyashenko and Yakovenko \cite{IYa}.

Our estimate depends on four parameters: $n$, $C$, $a_1$, $R$. Let
$F(x)=x^n+\sum\limits_{i=1}^{n-1}a_i x^i$, where $n=2l$ is the even
degree of the monic polynomial $F$ without a constant term,
$\forall i$ $|a_i|<C$, so $C$ is the size of a compact subset in
the space of parameters, $R$ is the size of the neighborhood of the
origin, such that there are no bigger than $l$ limit cycles located
outside of this neighborhood, $|a_1|$ stands the distance from the
equation linearization to the center case in the space of
parameters and $2-|a_1|$ stands the distance from the equation
linearization to the node case in the space of parameters.\newline
%\newline
\end{abstract}

\maketitle
\tableofcontents

\section{Hilbert-Smale problem}
In 1977 A.~Lins~Neto, W.~de~Melo and C.~C.~Pugh \cite{LMP}
examined small perturbations of a linear center for a special
class of polynomial vector fields on the plane. This class is
called Li\'enard equations:
\begin{equation}\label{eq:main}
\begin{cases}
\begin{aligned}
\dot x&{}=y-F(x),\\
\dot y&{}=-x,
\end{aligned}
\end{cases}
\end{equation}
where $F$ is a polynomial of odd degree. Actually, Li\'enard in
1928 introduced it for a modeling of the non-linear damping in
electric circuits \cite{L}. It was a generalization of the
famous Van der Pol equation \cite{V}.

Authors of \cite{LMP} proved the finiteness of limit cycles for a
Li\'enard equation of odd degree $n$. Let us remind that the
Finiteness problem (also known as the "Dulac problem") was solved
in full generality only in 1991 by Ilyashenko \cite{I1} and in
1992 by \'Ecalle \cite{E} independently.

Also A.~Lins~Neto, W.~de~Melo and C.~C.~Pugh \cite{LMP} conjectured
that the number of limit cycles of~$(\ref{eq:main})$ is not bigger
than $\frac{n-1}{2}$.

In 1998 S.~Smale \cite{S} suggested to consider a restriction of
the second part of the Hilbert's 16th problem to Li\'enard
equations of odd degree. He conjectured that there exists an
integer $n$ and real $C$ such that the number of limit cycles
of~$(\ref{eq:main})$ is not bigger than $Cn^q$.

In 1999 Yu.~Ilyashenko and A.~Panov \cite{IP} got an explicit
upper bound for the number of limit cycles of Li\'enard equations
through the (odd) degree of the monic polynomial $F$ and
magnitudes of its coefficients. Their result reclined on the
theorem of Ilyashenko and Yakovenko that binds the number of
zeros and the growth of a holomorphic function \cite{IYa}.

In 2007 F.~Dumortier, D.~Panazzolo and R.~Roussarie \cite{DPR}
constructed a counterexample to the conjecture of A.~Lins~Neto,
W.~de~Melo and C.~C.~Pugh. Namely, they presented an example of
a Li\'enard equation of odd degree $n$ with at least
$\frac{n+1}{2}$ limit cycles.

In 2008 Yu.~Ilyashenko \cite{I3} suggested to prove a result
analogous to the one of Ilyashenko and Panov for Li\'enard
equations of even degree.

In 2008 M.~Caubergh and F.~Dumortier in \cite{CD} proved the
following theorem for Li\'enard equations of even degree.

\begin{thm}\label{thm:CD}
Let $K$ be a compact set of polynomials of degree exactly $n=2l$,
then there exists $R>0$ such that any system having an
expression~$(\ref{eq:main})$ with $F\in K$ has at most $l$ limit
cycles having an intersection with $\Ro^2\backslash B_R$.
\end{thm}

Here and bellow $B_R$ denotes the ball around the origin with
the radius $R$.

\section{Notations and the Ilyashenko strategy}
From now on we will consider a system~$(\ref{eq:main})$, where $F$
is a monic polynomial of even degree $n=2l$ without a constant
term.

\begin{rem}
The assumption $F(0)=0$ does not reduce the generality; it may be
fulfilled by a shift $y\mapsto y+a$. The assumption that $F$ is
monic may be fulfilled by rescaling in $x$, $y$ and reversing the
time if necessary.
\end{rem}

Let $v$ be an analytic vector field in the real plane, that may be
extended to $\Co^2$. For any set $D$ in a metric space denote by
$U^{\eps}(D)$ the $\eps$-neighborhood of $D$. The metrics in $\Co$
and $\Co^2$ are given by:
\begin{align*}
&\rho(z,w)=|z-w|, &&z,w\in\Co;\\
&\rho(z,w)=\max(|z_1-w_1|,|z_2-w_2|), &&z,w\in\Co^2.
\end{align*}

Denote by $|D|$ the length of the segment $D$. For any larger
segment $D'\supset D$, let $\rho(D,\partial D')$ be the Hausdorff
distance between $D$ and $\partial D'$.

We want to apply the next theorem proved by Ilyashenko and
Panov \cite{IP}. In fact, it is the easy corollary from the
Growth-and-Zeros theorem for holomorphic functions proved by
Ilyashenko and Yakovenko \cite{IYa}.

Consider the system
\begin{equation}\label{eq:v.f.}
\dot{x}=v(x), \qquad \qquad x\in\Ro^2.
\end{equation}

\begin{thm}\label{thm:IP1}
Let $\Gamma$ be a cross-section of the vector field~$v$,
$D\subset\Gamma$ a segment. Let $P$ be the Poincar\'e map
of~$(\ref{eq:v.f.})$ defined on $D$, and ${D\subset D'=P(D)}$.
Suppose that $P$ may be analytically extended to\linebreak
${U=U^{\eps}(D)\subset\Co}$, $\eps<1$, and
${P(U)\subset U^1(D')\subset\Co}$. Then the number $\#LC(D)$ of
limit cycles that cross $D$ admits an upper estimate:
\begin{equation}\label{eq:IP-th}
\#LC(D)\leq
e^{2|D|\eps^{-1}}\log\frac{|D'|+2}{\rho(D,\partial D')}.
\end{equation}
The same is true for $P$ replaced by $P^{-1}$.
\end{thm}

Actually, the {\it Ilyashenko strategy} is the application of the
previous theorem. It requires a purely qualitative investigation
of a vector field, i.e. a construction of such $D$ for every nest
of limit cycles. This strategy was applied before in
papers~\cite{I2} and~\cite{IP}.

We take $K$ from the Theorem~$\ref{thm:CD}$ to be the space of
monic polynomials of degree exactly $n$ with coefficients, which
absolute values are bounded by some positive constant $C\geq4$,
i.e.
$$F(x)=x^n+\sum\limits_{i=1}^{n-1}a_i x^i, \qquad
\forall i: |a_i|<C.$$

If $0<|a_1|<2$ then the unique singular point $(0,0)$ of the
system~$(\ref{eq:main})$ is a focus. In our work we will consider
only this case.

% In fact, $R$ is a big constant. So, $R$ would be bigger than
% $C+1$. From now on we will replace $R$ by $\max(R,C+1)$.

\section{Bendixson trap from within}
\begin{figure}[ht]
  \includegraphics[width=7cm,height=7cm]
  {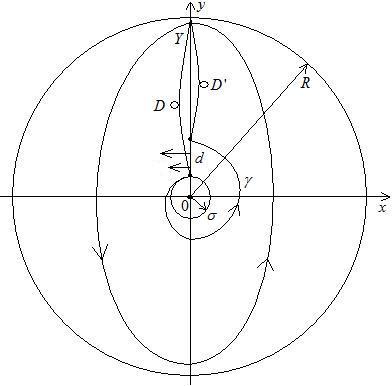}\\
  \caption{The inverse Poincar\'e map of the Li\'enard equation
  $(\ref{eq:main})$ inside the ball $B_{R}$.}
  \label{jpg:bendixson_from_within}
\end{figure}

In this Section we construct an interval $D$, which lies inside
$B_R$ and intersects transversally all limit cycles in
$B_R$. Also we find an upper estimate for the {\it Bernstein
index}, $b=\log\frac{|D'|+2}{\rho(D,\partial D')}$. To do that we
need to estimate $\rho(D,\partial D')$ from bellow, where
${D'=P(D)\subset D}$ and $P$ is the Poincar\'e map defined on $D$
(see the Figure~$\ref{jpg:bendixson_from_within}$).

Let $\varphi$, $r$ be polar coordinates on $\Ro^2$,
$\dot{\varphi}$, $\dot{r}$ be derivatives with respect
to~$(\ref{eq:main})$.

First of all we need to determine the size of the domain, there
the Poincar\'e map is defined.

\begin{lem}\label{lem:Poincare_domain}
Put $\sigma=\frac{|a_1|(2-|a_1|)}{8C}e^{\frac{8\pi}{|a_1|-2}}$.
In the focus case $\left(0<|a_1|<2\right)$ the Poincar\'e map for
the system~$(\ref{eq:main})$ is well defined in $B_{\sigma}$.
\end{lem}
\begin{proof}
Let us calculate $\dot{r}$.
\begin{multline}\label{eq:radius_dot}
\dot{r}=\frac{x\dot{x}+y\dot{y}}{r}=
\frac{r\cos\varphi(r\sin\varphi-F(r\cos\varphi))-
r^2\sin\varphi\cos\varphi}{r}=\\
=-\cos\varphi F(r\cos\varphi)=
-r\cos^2\varphi\sum_{i=1}^{n}a_i (r\cos\varphi)^{i-1}=\\
=-r\cos^2\varphi\left(a_1+O(r,\varphi)\right),
\end{multline}
where
$O(r,\varphi)=\sum\limits_{i=2}^{n}a_i (r\cos\varphi)^{i-1}$.

Let us calculate $\dot{\varphi}$.
\begin{multline}\label{eq:angle_dot}
\dot{\varphi}=\frac{x\dot{y}-y\dot{x}}{r^2}=
\frac{-r^2\cos\varphi-r^2\sin^2\varphi+
r\sin\varphi F(r\cos\varphi))}{r^2}=\\
=-1+\frac{\sin\varphi F(r\cos\varphi)}{r}=
-1+\sin2\varphi\left(\frac{a_1}{2}+\frac{O(r,\varphi)}{2}\right).
\end{multline}

The absolute value of the function $O(r,\varphi)$
admits the following upper estimate in $B_{\frac12}$:
\begin{equation}\label{eq:tail}
\left|O(r,\varphi)\right|\leq\sum_{i=2}^{n}Cr^{i-1}=
Cr\frac{1-r^{n-1}}{1-r}<\frac{Cr}{1-r}\leq2Cr.
\end{equation}

Therefore, in $B_{\frac{2-|a_1|}{4C}}\subset B_{\frac12}$:
$\dot{\varphi}\leq\frac{|a_1|-2}{4}$. Indeed,
\begin{multline*}
\dot{\varphi}\leq
-1+\left|\sin2\varphi\right|\left(\frac{|a_1|}{2}+
\frac{\left|O(r,\varphi)\right|}{2}\right)
\leq-1+\frac{|a_1|}{2}+Cr\leq\\
\leq\frac{|a_1|-2}{2}+C\frac{2-|a_1|}{4C}=\frac{|a_1|-2}{4}.
\end{multline*}

Also, in $B_{\frac{2-|a_1|}{4C}}$: $|\dot{r}|\leq2r$.
Indeed,
\begin{multline*}
|\dot{r}|\leq|r\cos^2\varphi|\left(|a_1|+|O(r,\varphi)|\right)
\leq\left(|a_1|+2Cr\right)r\leq\\
\leq\left(|a_1|+2C\frac{2-|a_1|}{4C}\right)r
\leq\frac{2+|a_1|}{2}r\leq2r.
\end{multline*}

Hence, any trajectory starting from any point from $B_{\sigma}$
rotates around the origin on the angle not less than $2\pi$ before
leaving $B_{\frac{2-|a_1|}{4C}}$.

Indeed, $\dot{\varphi}\leq\frac{|a_1|-2}{4}$ implies that during
the time, ${\vartriangle t=2\pi\frac{4}{2-|a_1|}}$ the variation of
the angle, ${\vartriangle\varphi\geq2\pi}$ and the variation of the
radius,\linebreak
${\vartriangle r\leq e^{2\frac{4\pi}{2-|a_1|}}}$ during the
same time ${\vartriangle t}$, because ${|\dot{r}|\leq2r}$.

Finally,
$\sigma e^{\frac{8\pi}{2-|a_1|}}=\frac{|a_1|(2-|a_1|)}{8C}
\leq\frac{2-|a_1|}{4C}$.
\end{proof}

Let us denote by $Y$ the maximal $y$-coordinate of the point of
intersection between the most external limit cycle which lies
inside $B_R$ (if it exists, of course) and $y$-axis.

\begin{lem}\label{lem:hausdorff}
If $a_1$ is negative, then $\dot{r}>0$ in $B_{\sigma}$. Let
$D=[\sigma,Y]\subset 0y$. Then
$d=\rho(D,\partial D')\geq\frac{\pi|a_1|}{2}\sigma$.
\end{lem}

\begin{proof}
If $r<\sigma$, then $r<\frac12$ and by~$(\ref{eq:tail})$:
$\left|O(r,\varphi)\right|\leq2Cr<\frac{|a_1|}{2}$.

Therefore by~$(\ref{eq:radius_dot})$,
$$\dot{r}>r\cos^2\varphi\left(-a_1-\frac{|a_1|}{2}\right)=
-\frac{a_1}{2}r\cos^2\varphi>0.$$

This proves the first part of the Lemma.

Consider the orbit $\gamma$ of the system~$(\ref{eq:main})$ that
passes through the point $(0,\sigma)$. Then the Hausdorff
distance, $d$ can be estimated as follows:
\begin{equation*}
d\geq\left|\int_0^{2\pi}\dot{r}(\gamma)d\varphi\right|>
\int_0^{2\pi}-\frac{a_1}{2}\sigma\cos^2\varphi d\varphi
=\frac{\pi|a_1|}{2}\sigma.
\end{equation*}

This inequality completes the proof of the Lemma.
\end{proof}

\begin{rem}
For positive $a_1$ we can get the same results just by reversing
of the time.
\end{rem}

Now we can estimate $b$ from above:
\begin{equation}\label{eq:bernstein}
b\leq\log\frac{R+2}{d}\leq\log\frac{2(R+2)}{\pi|a_1|\sigma}<
\frac{R+2}{|a_1|\sigma}.
\end{equation}

\section{Complex domain of the inverse Poincar\'e map}
The Theorem~$\ref{thm:IP1}$ uses the width $\eps$ of the complex
domain $U^{\eps}(D)$ to which the (inverse) Poincar\'e map may be
extended. We will apply the following theorem to estimate this
$\eps$ from bellow.

\begin{thm}\label{thm:eps}
Let $P:D\rightarrow D'$ be the Poincar\'e map of~$(\ref{eq:v.f.})$.
For any $x\in D$ denote by $\varphi_{x,P(x)}$ the arc of the phase
curve of $(\ref{eq:v.f.})$ starting at $x$ and ending at $P(x)$.

Let
$$\Omega(D)=\bigcup_{x\in D}\varphi_{x,P(x)},$$
and
\begin{equation}\label{eq:eps_mu_L}
1\leq\mu=\max_{U^2(\Omega)}|v|, \qquad\quad L=2\mu.
\end{equation}

Let $t(x)$ be the time length of the arc $\varphi_{x,P(x)}$, and
$$T_{\max}=\max_{\substack{x\in D}}t(x), \qquad T=T_{\max}+1.$$

Let
\begin{equation}\label{eq:eps_delta_lambda}
\delta\leq e^{-LT}, \qquad \lambda=\sqrt{\delta}, \qquad
\eps=\delta^2.
\end{equation}

Suppose that ${(z_1,z_2)}$ are coordinates in ${\Co^2}$,
${^{\Co}\Gamma=\{z_1=0\}}$,\linebreak
${v=(v_1,v_2)}$.

Let
${\Pi_{\delta}=U^{\delta}(0)\times U^{\lambda}(D')\subset\Co^2}$.
Suppose that
\begin{equation}\label{eq:eps_frac}
\left|\frac{v_2}{v_1}\right|\leq\mu
\quad \text{in} \quad \Pi_{\delta}.
\end{equation}

Then the Poincar\'e map ${P:D\rightarrow D'}$
of~$(\ref{eq:v.f.})$ may be analytically extended to
${U^{\eps}(D)\subset{^{\Co}\Gamma}}$, and
$P(U^{\eps}(D))\subset U^{1}(D)$.

The same is true for $P$ replaced by $P^{-1}$. In this case
$P^{-1}(D)=D'$,
$\Omega(D)=\bigcup\limits_{\substack{x\in D'}}\varphi_{x,P(x)}$.
\end{thm}
\begin{proof}[For the proof see~\cite{IP}]
\end{proof}

Bellow we will produce some preliminary calculations, which would
allow us to apply the Theorem~$\ref{thm:eps}$ later.

\begin{Def}\label{def:C-monic}
A {\it $C$-monic polynomial} is a real polynomial in one variable
with the highest coefficient one and other coefficients no greater
than $C$ in absolute value, with zero constant term.
\end{Def}

\begin{prop}[Properties of $C$-monic polynomials]
\label{prop:C-monic}
Let $F$ be a\linebreak
{$C$-monic} polynomial of degree $n$, $C\geq2$. Then
\begin{align}
&\max\limits_{x\in\left[0,X\right]}|F(x)|\leq2X^n &&\text{for }
X\geq C+1,\label{eq:C-monic_1}\\
&\max\limits_{x\in\left[0,X\right]}|F'(x)|\leq Cn^2X^{n-1}
&&\text{for } X\geq1,\label{eq:C-monic_2}\\
&|F(z)|\leq2C|z| &&\text{for } z\in\Co,
|z|\leq\frac12.\label{eq:C-monic_3}
\end{align}
\end{prop}
\begin{proof}[For the proof see~\cite{IP}]
\end{proof}

\begin{lem}\label{lem:mu_L}
Let $v$ be the vector field given by the system~$(\ref{eq:main})$.
Then $\mu$ and $L$ from the Theorem~$\ref{thm:eps}$ admits the
following estimates:
\begin{equation}\label{eq:est_mu_L}
\mu\leq3(R+2)^n \qquad \qquad L\leq6(R+2)^n.
\end{equation}
\end{lem}
\begin{proof}
By definition, $U^2(\Omega)\subset B_{R+2}$. So
$$|v|\leq|\dot{x}|+|\dot{y}|\leq|x|+|y|+|F(x)|\leq
2(R+2)+2(R+2)^n,$$
where the last inequality provided by~$(\ref{eq:C-monic_1})$.
Hence,
$$\mu\leq2\left(R+2+(R+2)^n\right)\leq3(R+2)^n,
\qquad L=2\mu\leq6(R+2)^n,$$
that proves the Lemma.
\end{proof}

Let $G=B_R\setminus B_{\sigma}$. Then
$\Omega=\bigcup\limits_{x\in D}\varphi_{x,P(x)}\subset G$.

\begin{lem}\label{lem:T_max}
Let $\gamma_y$ be the arc $\varphi_{y,P{-1}(y)}$ of the phase
curve of~$(\ref{eq:main})$, where $y\in D$. Then $t(y)$, the time
length of $\gamma_y$, admits an estimate
\begin{equation}\label{eq:est_T_max}
T_{\max}=\max\limits_{y\in D}t(y)\leq\frac{25C^2n^2R^n}{\sigma}.
\end{equation}
\end{lem}
\begin{proof}
The arcs $\gamma_y$, $y\in D$ belongs to $G$. We will split $G$
into two domains: $|\dot{x}|\leq\alpha$ and $|\dot{x}|>\alpha$ for
$\alpha$ small to be chosen later. The second domain contains two
parts of $\gamma_y$: one with $\dot{x}<-\alpha$, the other with
$\dot{x}>\alpha$. The time length of any of them is no greater
than $\frac{2R}{\alpha}$. In the next Proposition we will choose
$\alpha$ so small that the curvilinear strip
$$S_{\alpha}=\left\{(x,y)\in G:|y-F(x)|\leq\alpha\right\}$$
is crossed by the orbits of~$(\ref{eq:main})$ in the time no
greater than $1$.

\begin{prop}\label{prop:alpha_omega}
Let
\begin{equation}\label{eq:alpha_omega}
\omega=\frac{\sigma}{3C}, \qquad
\alpha=\frac{\omega}{2Cn^2R^{n-1}}=\frac{\sigma}{6C^2n^2R^{n-1}}.
\end{equation}
Then the time length of any arc of the orbit of~$(\ref{eq:main})$
located in $S_{\alpha}$ is no greater than $1$.
\end{prop}
\begin{proof}
By the symmetry arguments it is sufficient to prove that in
${S_{\alpha}^+=S_{\alpha}\cap\{x>0\}}$:
$$\frac{d}{dt}(y-F(x))\leq-2\alpha.$$
Let us first prove that in $S_{\alpha}$ we have: $|x|>\omega$.
Namely, let $|x|\leq\omega$, $|y-F(x)|\leq\alpha$. Then
$(x,y)\in D_{\sigma}$. Indeed,
$$|x|+|y|\leq\omega+\alpha+
\max\limits_{\left[0,\omega\right]}|F(x)|.$$
By~$(\ref{eq:alpha_omega})$, $\alpha<\omega<\frac12$.
By~$(\ref{eq:C-monic_3})$, $|F(x)|\leq2C\omega$. Hence, for
$x\in\left[0,\omega\right]$,
$$|x|+|y|\leq(2C+2)\omega<3C\omega=\sigma.$$
By~$(\ref{eq:C-monic_2})$, $|F'(x)|\leq Cn^2R^{n-1}$ in $G$.
Therefore, for $x$ such that ${(x,y)\in S_{\alpha}^+}$ we have:
$x>\omega$, and
$$\frac{d}{dt}(y-F(x))=-x-F'(x)(y-F(x))
\leq-\omega+\alpha Cn^2R^{n-1}<-2\alpha,$$
because $\qquad$ $\qquad$
$\alpha=\dfrac{\omega}{2Cn^2R^{n-1}}<
\dfrac{\omega}{2Cn^2R^{n-1}+2}$.
\end{proof}

Let us finish the proof of the Lemma~$\ref{lem:T_max}$.

The arc $\gamma_y$ spends in $S_{\alpha}$ no longer time than $2$
(two crossings, each one no longer in time than 1, by the previous
Proposition); in $G\setminus S_{\alpha}$ no longer time than
$\frac{4R}{\alpha}$ (two crossings, one to the left, another to
the right with $|\dot{x}|\geq\alpha$). Hence,
$$T_{\max}\leq2+\frac{4R\cdot6C^2n^2R^{n-1}}{\sigma}<
\frac{25C^2n^2R^n}{\sigma}.$$

This calculation completes the proof of the
Lemma~$\ref{lem:T_max}$.
\end{proof}

\begin{rem}\label{rem:est_T}
The same inequality holds for $T_{\max}$ replaced by
${T_{\max}+1}$.
\end{rem}

Let us check the last assumption of the Theorem~$\ref{thm:eps}$.
\begin{lem}\label{lem:est_v.f.}
Take
\begin{equation}\label{eq:eps_prec}
\eps=\exp\left(-\frac{300C^2n^2R^n(R+2)^n}{\sigma}\right),
\quad \delta=\sqrt{\eps}, \quad \lambda=\sqrt{\delta}.
\end{equation}
Let, as in the Theorem~$\ref{thm:eps}$,
${\Pi_{\delta}=U^{\delta}(0)\times U^{\lambda}(D')\subset\Co^2}$.
Then in $\Pi_{\delta}$:
$$\left|\frac{v_2}{v_1}\right|<\mu.$$
\end{lem}
\begin{proof}
By~$(\ref{eq:C-monic_3})$ and by definition of $\Pi_{\delta}$,
$$|v_1(z)|\geq\left||z_1|-F(z)\right|\geq
(\sigma-\lambda)-2C\delta>\delta,$$
where the last inequality is trivial.
On the other hand, $v_2=-x$. In $\Pi_{\delta}$,
$|v_2|\leq\delta$. Hence, $\left|\frac{v_2}{v_1}\right|<1<\mu.$
\end{proof}

\begin{lem}\label{lem:eps_verif}
The inverse Poincar\'e map of the Li\'enard
equation~$(\ref{eq:main})$ may be extended to the domain
$U^{\eps}(D)\subset\Co$, where
$\eps=e^{-\frac{300C^2n^2R^n(R+2)^n}{\sigma}}$.
Moreover, $P^{-1}(U^{\eps}(D))\subset U^{1}(D')$.
\end{lem}
\begin{proof}
This Lemma following from the Theorem~$\ref{thm:eps}$.
Lemmas~$\ref{lem:mu_L}$ and $\ref{lem:est_v.f.}$ verifies
assumptions~$\ref{eq:eps_mu_L}$ and $\ref{eq:eps_frac}$
respectively. We only should check the
assumption~$\ref{eq:eps_delta_lambda}$. By the
Remark~$\ref{rem:est_T}$, $T<\frac{25C^2n^2R^n}{\sigma}.$
Hence,
\begin{multline*}
\delta=\sqrt{\eps}=
\exp\left(-\frac{150C^2n^2R^n(R+2)^n}{\sigma}\right)=\\
=\exp\left(-6(R+2)^n\frac{25C^2n^2R^n}{\sigma}\right)\leq e^{-LT},
\end{multline*}
that proves the Lemma.
\end{proof}

\section{Final estimate}
\begin{thm}
The number $L(n,C,a_1,R)$ of limit cycles of~$(\ref{eq:main})$ in
the case when $n$ is even, $C\geq4$ and $0<|a_1|<2$, admits the
following upper bound:
\begin{equation*}
L(n,C,a_1,R)<
\exp\left(\exp\left(\frac{38400C^4n^2R^{n+1}(R+2)^{n+1}}
{|a_1|^3(2-|a_1|)^2}e^{\frac{16\pi}{2-|a_1|}}\right)\right).
\end{equation*}
\end{thm}

\begin{proof}
Now we can apply the Theorem~$\ref{thm:IP1}$.
By definition, $|D|$ and $|D'|$ are less than $R$. The
Lemma~$\ref{lem:eps_verif}$ provides us with the lower bound on
$\eps$. So estimates~$(\ref{eq:IP-th})$
and~$(\ref{eq:bernstein})$ imply:
\begin{multline*}
L(n,C,a_1,R)<
\exp\left(2R\exp\left(\frac{300C^2n^2R^n(R+2)^n}
{\sigma}\right)\right)\frac{R+2}{|a_1|\sigma}<\\
<\exp\left(\frac{2R(R+2)}{|a_1|\sigma}
\exp\left(\frac{300C^2n^2R^n(R+2)^n}{\sigma}\right)\right)<\\
<\exp\left(\exp\left(\frac{600C^2n^2R^{n+1}(R+2)^{n+1}}
{|a_1|\sigma^2}\right)\right)=\\
=\exp\left(\exp\left(\frac{38400C^4n^2R^{n+1}(R+2)^{n+1}}
{|a_1|^3(2-|a_1|)^2}
\exp\left(\frac{16\pi}{2-|a_1|}\right)\right)\right).
\end{multline*}
This calculation completes the proof of the Theorem.
\end{proof}

\section{Acknowledgments}
The author is grateful to D.~V.~Anosov for the scientific
advising, to Yu.~S.~Ilyashenko for posing the problem and for
the constant attention to this work and to the Fields Institute
(Toronto, Ontario, Canada) for the hospitality and excellent
working conditions.

\end{document}